# Novel closed-loop controllers for fractional linear quadratic tracking systems

Iman Malmir[†]

**Abstract** A new method for finding closed-loop optimal controllers of fractional tracking quadratic optimal control problems is introduced. The optimality conditions for the fractional optimal control problem are obtained. Illustrative examples are presented to show the applicability and capabilities of the method.

## 1 Introduction

Fractional optimal tracking control as a combination of fractional optimal control and tracking control [1,2], aims at finding the optimal control law to make the fractional system tracks the desired reference in an optimal way. The tracking system is widely used in aerospace and mechanical systems. For example, consider an antenna control system to track an aircraft. In a huge system, the complexity of optimization problems increases exponentially [3] and the traditional methods are not applicable. In [4], we used a new idea together with orthogonal functions to find the closed loop solutions of fractional linear quadratic regulator problems. Orthogonal functions are fundamental in solving a wide variety of engineering and mathematical problems, including differential and integral equations, approximation theory, and dynamical systems [5]. In this paper, a feedback optimal control for the fractional optimal control problem that involves a quadratic Bolza tracking performance index and a fractional linear state equation is proposed.

Fractional models describe the behavior of some systems better than the integer systems [6]. Fractional model is now a well-established tool in engineering science [7]. Fractional differential equations are applied in Engineering, science, finance and bioengineering [8]. In engineering systems, we are dealing with tracking problems, for example, see [9,10].

[†]Aerospace and Mechanical Engineering Group, Ronin Institute, Montclair, USA. E-mail: iman.malmir@outlook.com

There are two different perspectives in designing the feedback controls of the fractional systems. In the first, we design the controller, then by adjusting some parameters we find the elements of this controller. The first perspective has been studied in many papers. In the second perspective, we are dealing with fractional optimal control of a specified performance index based on objectives in which the feedback controller must satisfy some conditions and equations; it should be noted that only few works have studied the closed-loop solutions of the fractional optimal control problems as the second perspective. There are also too many texts in the fractional optimal control that have only researched open-loop solutions. For more details, see [11].

Many significant real world challenges arise as optimization problems on different classes of control systems [12]. For example, computational models to uncover the physiological mechanisms underlying the disease-related neural circuits receiving Parkinson's disease was developed in [13]. The closed loop controller for fractional tracking linear quadratic was presented in a few texts.

In [6], the explicit formula of the optimal tracking control-state pair of tracking control problem for fractional order descriptor systems was presented. An application of the fractional tracking optimal control, in controlling the glucose levels in the blood for diabetics was presented in [14]. Here, we are going to use the results given in [4,15,16] for finding the closed-loop optimal controllers of fractional linear tracking quadratic optimal control problems. By using the results of [11,17], we can expand our idea to nonlinear tracking systems.

In recent texts, the solutions of the fractional optimal control problems were investigated. In [18], the Pontryagin maximum principle was obtained for the fractional optimal control problem with terminal and running state constraints and the necessary optimality conditions for the problem were derived. In [19], singular linear Volterra integral control system with quadratic cost functional over the infinite time horizon was studied and existence of the optimal control was discussed.

Fractional optimal control problems are generally defined by changing the integer derivative to a fractional derivative of desired senses [20–22]. Different classes of the fractional optimal control problems were studied in [23–25]. We saw in defining new fractional systems that there may be inconsistencies, for instance, see Sec 4.2.1 in [24]. Solving fractional optimal control problems is certainly more complicated than the old integer-order problems; in this research, we implement the method on fractional problems despite this complexity.

There are theoretical works that only study optimal control problems theoretically without solving even a simple example; but in reality, one can see the fact that it is impossible to find the optimal control law for fractional systems by using any of those methods. In the optimal control theory, each solution method should be tested and executed on real optimal control problems, for example, see [1] and [2]. Our method is presented theoretically and practically. There are many texts that only deal with solutions of integer optimal control problems theoretically in which anyone can see the fact that even a simple first-order optimization problem has not been solved to verify their methods.

It is practically useless to present pure theory without application because optimal control problems as models of the real-world systems need practical and efficient solution methods. Such works only study optimal control problems theoretically and give many formulas without solving a simple optimization problem, for example, see [26], but, in practice, their methods do not even have the ability to solve the integer order (the old version) of a typical optimal control problem like Example 6 in [27], let alone difficult problems like those in [24].

## 2 Problems statements

**Definition 1.** *The left Caputo derivative of order $\alpha$ for a function $f(t)$ "${}^C_0D^\alpha_t f(t)$" is defined as*

$$
{}^C_0D^\alpha_t f(t) := \begin{cases} \frac{1}{\Gamma(n-\alpha)} \int_0^t (t-\rho)^{n-\alpha-1} f^{(n)}(\rho) \, d\rho, & n-1 < \alpha < n \\ \frac{d^n}{dt^n} f(t), & \alpha = n, \end{cases}
$$

*where $n \in \mathbb{N}$.*

**Remark 1.** *The left Caputo derivative of order $\alpha \in (0,1[$ takes the form*

$$
{}^C_0D^\alpha_t f(t) := \frac{1}{\Gamma(1-\alpha)} \int_0^t (t-\rho)^{-\alpha} \dot{f}(\rho) \, d\rho.
$$

*Obviously, one can see that the term $(t-\rho)^{-\alpha}$ exists inherently in the integral.*

### 2.1 Fractional quadratic linear tracking problem

The fractional linear tracking problem consists of the plant described by

$$
{}^C_0D^\alpha_t \mathbf{x}(t) = \mathbf{A}(t)\mathbf{x}(t) + \mathbf{B}(t)\mathbf{u}(t), \quad 0 \leq t \leq t_f \tag{1}
$$

the real-valued performance index

$$
J = \tfrac{1}{2}[\mathbf{x}(t_f) - \mathbf{r}(t_f)]^\top \mathbf{T}[\mathbf{x}(t_f) - \mathbf{r}(t_f)] + \tfrac{1}{2}\int_0^{t_f} \{[\mathbf{x}(t) - \mathbf{r}(t)]^\top \mathbf{Q}[\mathbf{x}(t) - \mathbf{r}(t)] + \mathbf{u}^\top(t)\mathbf{R}\mathbf{u}(t)\} \, dt, \tag{2}
$$

and the initial condition as $\mathbf{x}(0) = \mathbf{x}_0$, where the final time $t_f$ is fixed, and $\mathbf{x}(t_f)$ is free. Here, ${}^C_0D^\alpha_t$ is the left Caputo derivative of order $\alpha$, and $\alpha \in (0,1]$, $\mathbf{x}: [0, t_f] \to \mathbb{R}^q$ and $\mathbf{u}: [0, t_f] \to \mathbb{R}^r$ are the state and control vectors, $\mathbf{A}(t) \in \mathbb{R}^{q \times q}$, $\mathbf{B}(t) \in \mathbb{R}^{q \times r}$, $F: \mathbb{R} \times \mathbb{R}^q \times \mathbb{R}^r \to \mathbb{R}$ is a real-valued function, $\mathbf{r}: [0, t_f] \to \mathbb{R}^q$ is the desired or reference state vector and is continuous on its interval, $\mathbf{T}, \mathbf{Q}(t) \in \mathbb{R}^{q \times q}$ are symmetric, positive semi-definite matrices, and $\mathbf{R}(t) \in \mathbb{R}^{r \times r}$ a is symmetric, positive definite matrix. $\mathbf{A}(t)$, $\mathbf{B}(t)$, $\mathbf{Q}(t)$, and $\mathbf{R}(t)$ can be either time-varying or time-invariant. Our objective is to control this system in such a way that the state $\mathbf{x}(t)$ tracks the desired state $\mathbf{r}(t)$ as close as possible during the time interval $[0, t_f]$ with minimum expenditure of control effort. In other words, the optimal tracking problem is to find the closed-loop optimal control $\mathbf{u}^*(t)$ for the fractional system (1) such that the performance index in (2) is minimized.

**Remark 2.** *In tracking optimal control problems, we are interested to minimize the error between the actual state and the reference or the desired state, hence we consider the square term '$(\mathbf{x}(t) - \mathbf{r}(t))^2$' in the performance index to minimize the tracking error. Although we wish to the optimal state be equal to the desired reference, but, in reality, we have*

$$\mathbf{x}^*(t) \neq \mathbf{r}(t).$$

*For achieving the goal, we design the values of the weighting matrices, for example, see [9]. In addition, if one takes $\mathbf{x}^*(t) - \mathbf{r}(t) \equiv \mathbf{0}$, considering 0 in the integrand of the performance index* (2) *does not provide any realistic sense.*

## 2.2 Fractional quadratic nonlinear tracking problem

The fractional nonlinear tracking problem consists of the plant described by

$$\,^C_0 D^\alpha_t \mathbf{x}(t) = \mathbf{f}(\mathbf{x}(t), \mathbf{u}(t), t), \tag{3}$$

and the performance index or the cost functional is presented in (2).

## 2.3 Pure fractional quadratic linear/nonlinear tracking problem: open problem

In some works, the performance index of the fractional optimal control problem is also defined by a fractional order, for example, [23–25]. As an open problem, we can consider the following performance index for the fractional linear or nonlinear tracking system (1) and (3), instead of (2),

$$\begin{aligned} J(\mathbf{u}(t), \alpha_1) = &\tfrac{1}{2} [\mathbf{x}(t_f) - \mathbf{r}(t_f)]^\top \mathbf{T} [\mathbf{x}(t_f) - \mathbf{r}(t_f)] \\ &+ \tfrac{1}{2} \,^{RL}_{0} I^{\alpha_1}_{t_f} \left\{ [\mathbf{x}(t) - \mathbf{r}(t)]^\top \mathbf{Q} [\mathbf{x}(t) - \mathbf{r}(t)] + \mathbf{u}^\top(t) \mathbf{R} \mathbf{u}(t) \right\}, \end{aligned} \tag{4}$$

where, the constant Riemann–Liouville integral of order $\alpha$ for a function $f(t)$ denoted by "$^{RL}_{t_0} I^\alpha_{t_f} f(t)$" in (4) is defined for $t_0 \leq t \leq t_f$ as [24]

$$^{RL}_{t_0} I^\alpha_{t_f} f(t) = \frac{1}{\Gamma(\alpha)} \int_{t_0}^{t_f} (t_f - t)^{\alpha-1} f(t)\, dt. \tag{5}$$

We start with the results presented in [4] for the fractional regulator systems. Then, by using the transformation and the procedure given in [16], and using the results presented in [15], we obtain the optimality conditions for the first problem. After that, by using the results presented in [11] and the mentioned transformation, we can solve the second fractional optimal control problem.

**Remark 3.** *In the optimal control theory, we describe a plant by a set of linear differential equations as a linear regulator problem. The system can be a liquid-level control system, a mass-spring-damper system, an automobile suspension system, an n-plate gas absorber system, a stirred tank chemical reactor, a system of treatment for a disease and etc. Hence, one can see the fact that* (1) *is the state equation of* *the control systems and its state or control may represent the temperature, pressure, density, deviation from concentration of a system, weight fractions of species and inlet feed rate of components of a tank reactor, the drug concentration of a disease model and so on.*

## 3 Closed-loop optimal controls and optimality conditions for the fractional tracking problems

**Lemma 1.** *If functions $f(t)$ and $g(t)$ are in $C[0, t_f]$, then for $\alpha \in (0.9, 1]$*

$$\int_0^{t_f} f(t) \, {}_0^C D_t^\alpha g(t) \, dt \simeq \Delta_{0f} - \int_0^{t_f} ({}_0^C D_t^\alpha f(t)) g(t) \, dt, \tag{6}$$

*where $\Delta_{0f}$ denotes constant terms.*

*Proof.* Taking $\Delta_{0f} = \Delta_{0f}^0 + \Delta_{0f}^1 + \Delta_{0f}^2$, the proof was presented in [4]. □

The following points are worth noting:

- For $\alpha = 1$, the proof is easy and we did not present it. Also, since for $\alpha = 1$, we have '=' in (6) and for $\alpha \in (0.9, 1)$, we used '≈' in the approximation after ninth equality in the proof of Lemma 3.1 in [4], we utilized '≃' in (6).

- We used the calculus of the variations to find the optimality conditions in [4] as a common approach to find these conditions in the optimal control theory. There are constant terms in the given lemma, but they are used for boundary conditions and we do not need them in our method. One can refer to Pages 230–233 in [1] and Pages 61–65 in [2] as perfect sources in the optimal control theory, to see how boundary conditions are derived, for example, see Eq. (2.7.32) on Page 65 of [2]. The method in [4] is a hybrid method (combination of two methods, one for **K** and for **l**, and one for removing boundary conditions) and boundary conditions are not derived because, by using another method and combining the approaches, we do not need them. The terms have no impact on the optimality conditions, we did not discuss them in our approach.

- For the term $f(t) \frac{1}{\Gamma(1-\alpha)} \left[ (t-\rho)^{-\alpha} g(\rho) \right] \big|^{\rho=t}$, by using the expansions given in [4] for $\alpha \in (0.9, 1)$, we can write

$$\frac{1}{\Gamma(1-\alpha)} f(t) \left[ (t-\rho)^{-\alpha} g(\rho) \right] \big|^{\rho=t} \approx \frac{1}{\Gamma(1-\alpha)} f(t) \left[ (t-\rho)^{1-\alpha} g(\rho) \sum_{m=0} \frac{\rho^m}{t^{m+1}} \right] \big|^{\rho=t}$$

or

$$\approx \frac{1}{\Gamma(1-\alpha)} \Big[ \sum_{m_1=0} \sum_{m_2=0} \sum_{m_3=0} (t-\rho)^{1-\alpha} f_{m_1} g_{m_2} t^{m_1-m_3-1} \rho^{m_2+m_3} \Big] \big|^{\rho=t}$$

and since the term is constant and is used for finding boundary conditions, we did not analyze it. In other relations, it is obvious from Remark 1 that the term $(t - \rho)^{-\alpha}$ for $\rho \to t$ exists inherently in the integrals.

- In the analysis of $\delta J_a(\mathbf{u}^*(t)) = 0$ given in [4], the constants in the integrals are definite for 0 and $t_f$, they are constant values and have no effect on the optimality conditions, regardless of whether they depend on the functions $f$ and $g$. The definite integrals are constants and we have not discussed them because they have no effect on the optimality conditions given in [4] as Eq. (17). One can easily see in our method that there is no need to derive the boundary conditions for the state and costate variables.

- Such conditions (initial and final conditions) are used in the traditional method to find the optimal state and costate from their differential equations (for example, see Table 2.1, Step 4 in [2]), which are unnecessary in the method of [4]. Hence, we did not discuss the constant terms in [4].

## 3.1 Optimality conditions for fractional quadratic linear tracking problem

**Theorem 1.** *The optimality conditions for fractional quadratic linear tracking problem given in* (1) *and* (2)*, and for* $\alpha \in (0.9, 1]$ *are*

$$\mathbf{Q}(t)(\mathbf{x}^*(t) - \mathbf{r}(t)) + \mathbf{A}^\top(t)\boldsymbol{\lambda}^*(t) + {}_0^C D_t^\alpha \boldsymbol{\lambda}^*(t) \simeq \mathbf{0}, \qquad (7)$$

$$\mathbf{R}(t)\mathbf{u}^*(t) + \mathbf{B}^\top(t)\boldsymbol{\lambda}^*(t) = \mathbf{0} \qquad (8)$$

*and*

$$\mathbf{A}(t)\mathbf{x}^*(t) + \mathbf{B}(t)\mathbf{u}^*(t) - {}_0^C D_t^\alpha \mathbf{x}^*(t) = \mathbf{0}, \qquad (9)$$

*where for the Riccati matrix* $\mathbf{P}(t)$ *and a vector* $\mathbf{z}(t)$*, we take* $\boldsymbol{\lambda}^*(t) = \mathbf{P}(t)\mathbf{x}^*(t) - \mathbf{P}(t)\mathbf{r}(t) + \mathbf{z}(t)$. $\mathbf{P}(t)$ *is the solution of the Riccati fractional order differential matrix equation as*

$$\begin{aligned}{}_0^C D_t^1(\mathbf{P}(t))\mathbf{x}^*(t) =& {}_0^C D_t^{1-\alpha}[-\mathbf{Q}(t)\mathbf{x}^*(t) - \mathbf{A}^\top(t)\mathbf{P}(t)\mathbf{x}^*(t)] \\ & - \mathbf{P}(t) {}_0^C D_t^{1-\alpha}[\mathbf{A}(t)\mathbf{x}^*(t) - \mathbf{B}(t)\mathbf{R}^{-1}(t)\mathbf{B}^\top(t)\mathbf{P}(t)\mathbf{x}^*(t)], \end{aligned} \qquad (10)$$

*and for co-reference* $\mathbf{v}(t)$*,* $\mathbf{z}(t)$ *is the solution of the vector fractional differential equation*

$${}_0^C D_t^1 \mathbf{z}(t) = -{}_0^C D_t^{1-\alpha}[\mathbf{A}^\top(t)\mathbf{z}(t)] + \mathbf{P}(t){}_0^C D_t^{1-\alpha}[\mathbf{B}(t)\mathbf{R}^{-1}(t)\mathbf{B}^\top(t)\mathbf{z}(t) - \mathbf{v}(t)]. \qquad (11)$$

*Since for* $\alpha = 1$*, we have* '=' *in* (7) *and for* $\alpha \in (0.9, 1)$*, we have* '$\approx$'*, we use* '$\simeq$' *in* (7)*.*

*Proof.* First, we use the transformation given in [16]. Then, by applying it, we remodel the system (1). Now we define the co-reference vector $\mathbf{v}(t)$ where $\mathbf{A}(t)\mathbf{r}(t)$ is known and by using a series approximation we can find ${}_0^C D_t^\alpha \mathbf{r}(t)$ from Mittag–Leffler function, for example, see [19]. Thus, the new system is

$${}_0^C D_t^\alpha \bar{\mathbf{x}}(t) = \mathbf{A}(t)\bar{\mathbf{x}}(t) + \mathbf{B}(t)\mathbf{u}(t) + \mathbf{v}(t). \qquad (12)$$

By using the discussions in [4] and [15], we can reach the optimality conditions and the fractional differential equations given in (7)–(9) and (10), (11). □

**Lemma 2.** *Assume that in* (12) *for constants* $\mu_1, \mu_2, \mu_3 > 0$*,* $\|\mathbf{A}(t)\| \leq \mu_1$*,* $\|\mathbf{v}(t)\| \leq \mu_3$ *and by* $\mathbf{C}(t) := -\mathbf{B}(t)\mathbf{R}^{-1}(t)\mathbf{B}^\top(t)\mathbf{P}(t)$*,* $\|\mathbf{C}(\bar{\mathbf{x}}(t), t) - \mathbf{C}(\hat{\mathbf{x}}(t), t)\| \leq \mu_2 \|\bar{\mathbf{x}} - \hat{\mathbf{x}}(t)\|$*, where* $\bar{\mathbf{x}}, \hat{\mathbf{x}} \in \mathbb{R}^q$*. Then* (12) *has a unique solution.*

*Proof.* The proof is similar to that presented for Lemma 3.5 in [4]. □

## 3.2 Closed-loop optimal control for fractional quadratic linear tracking problem

The optimal control law for fractional quadratic linear tracking problem is given by

$$\mathbf{u}^*(t) = -\mathbf{K}(t)\mathbf{x}^*(t) + \mathbf{l}(t), \tag{13}$$

where

$$\mathbf{K}(t) = \mathbf{R}^{-1}(t)\mathbf{B}^\top(t)\mathbf{P}(t),$$

$\mathbf{K}(t)$ is the Kalman gain and the Riccati matrix $\mathbf{P}(t)$ is the solution of the matrix fractional differential Riccati equation given in (10). In traditional methods, the term $\mathbf{l}(t)$ in (13) should be derived from the given fractional differential equation (11). But, we use the hybrid method presented in [15] to find this term without solving the mentioned equation and in addition, we present it in a new form.

## 4 Simulation studies

In this section, we apply the proposed method for finding the closed-loop optimal control of the fractional tracking optimal control problem.

### 4.1 Example 1

Consider the fractional tracking Van der Pol oscillator problem, which is constructed from Van der Pol oscillator equation [28]. The system is given by

$$^C_0D^\alpha_t x_1(t) = x_2(t),$$
$$^C_0D^\alpha_t x_2(t) = -x_1(t) + (1 - x_1^2(t))x_2(t) + u(t),$$

where $x_1(0) = 1$ and $x_2(0) = 0$. The performance index for tracking problem is

$$J = \tfrac{1}{2}\int_0^5 \left\{(x_1(t) - r(t))^2 + x_2^2(t) + u^2(t)\right\}dt,$$

where $r(t) = (1 - 0.4t)\cos(t)$. The problem is finding the optimal control law.

**Remark 4.** *This problem was not originally for tracking purposes, but we modify it to implement the proposed method. One should note that for better tracking, we have to increase the values of the weighting matrix $\mathbf{Q}$, for example, see Examples 1–3 in [16], and $\mathbf{Q}$ in the examples given in [29].*

We define the optimal control law by

$$u^*(t) = -[k_{11}(t),\ k_{12}(t)]\mathbf{x}^*(t) + [k^r_{11}(t),\ 0][1 + (1 - 0.4t)\cos(t); 0]$$

and the optimal states, the optimal control and the feedback gains are shown in Figure 1 for $\alpha = 0.9$. Considering Remark 4, we change the value $q_{11}$ of the state weighting matrix and set $10(x_1(t) - r(t))^2$ in the performance index, that is,

$$J = \tfrac{1}{2}\int_0^5 \left\{10(x_1(t) - r(t))^2 + x_2^2(t) + u^2(t)\right\}dt.$$

The results are shown in Figure 2 and, as can be seen, by this new value of the state weighting matrix, the system can better track the desired reference $r(t)$.

## 4.2 Example 2

This problem is a fractional tracking version of the example was studied in [4, 30]. The fractional state equation is

$${}^C_0D^\alpha_t \mathbf{x}(t) = \begin{bmatrix} \mathbf{0} & \mathbf{I} \\ -\mathbf{M}^{-1}\boldsymbol{\kappa} & \mathbf{0} \end{bmatrix} \mathbf{x}(t) + \begin{bmatrix} 0 & 0 & \cdots & 0 & 1/m_L \end{bmatrix}^\top u(t),$$

The performance index is given as

$$J = \int_0^{10} \left( [\mathbf{x}(t) - \mathbf{r}(t)]^\top \begin{bmatrix} \boldsymbol{\kappa} & \mathbf{0} \\ \mathbf{0} & \mathbf{M} \end{bmatrix} [\mathbf{x} - \mathbf{r}(t)] + u^2(t) \right) dt,$$

where

$$\mathbf{r}(t) = [0; 0; 0; -\tfrac{13}{450}t^2 + \tfrac{163}{450}t - \tfrac{1}{3}; 0; 0; 0; 0; 0; 0],$$

$$\mathbf{M} = \begin{bmatrix} m_1 & 0 & \cdots & 0 \\ 0 & m_2 & \cdots & 0 \\ \vdots & \vdots & \ddots & \vdots \\ 0 & 0 & \cdots & m_L \end{bmatrix},$$

$$\boldsymbol{\kappa} = \begin{bmatrix} \kappa_1 + \kappa_2 & -\kappa_2 & 0 & 0 & \cdots & 0 & 0 & 0 \\ -\kappa_2 & \kappa_2 + \kappa_3 & -\kappa_3 & 0 & \cdots & 0 & 0 & 0 \\ 0 & -\kappa_3 & \kappa_3 + \kappa_4 & -\kappa_4 & \cdots & 0 & 0 & 0 \\ \vdots & \vdots & \vdots & \vdots & \ddots & \kappa_{L-2} + \kappa_{L-1} & -\kappa_{L-1} & 0 \\ 0 & 0 & 0 & 0 & \cdots & -\kappa_{L-1} & \kappa_{L-1} + \kappa_L & -\kappa_L \\ 0 & 0 & 0 & 0 & \cdots & 0 & -\kappa_L & \kappa_L \end{bmatrix}.$$

The initial conditions are

$$x_L(0) = 1, x_l(0) = 0, l = 1, 2, \ldots, L-1, L+1, \ldots, 2L.$$

The problem is to find the optimal control law for the fractional tracking system. We take $\alpha = 0.95$, $L = 5$, $m_l = 10$ kg, and $\kappa_l = 1$ N/m for $l = 1, 2, \ldots, L$. Hence, $\mathbf{x}(0) = [0; 0; 0; 0; 1; 0; 0; 0; 0; 0]$ m.

**Remark 5.** *This system was not originally for tracking purposes. One should note that for better tracking of any tracking optimal control problem, we have to try various values of its weighting matrices, for example, see the values of weighting matrices given in the numerical experiments of [16] and [29].*

For this fractional tracking system, by taking the optimal control law as

$$u^*(t) = -\mathbf{K}(t)\mathbf{x}^*(t) + \mathbf{K}^r(t)[0; 0; 0; -\tfrac{13}{450}t^2 + \tfrac{163}{450}t - \tfrac{4}{3}; 0; 0; 0; 0; 0; 0],$$

the optimal states and control, and the feedback gains are shown in Figure 3.

The closed-loop solutions of these problems cannot be obtained by many of available theoretical methods.

# References

[1] D.E. Kirk, Optimal Control Theory: An Introduction, Courier Corporation, 2004.

[2] D.S. Naidu, Optimal control systems, Idaho State University, USA, CRC Press, 2003.

[3] Y. Chen, Y. Nan, X. Sun, T. Tan, Superior Control of Spacecraft Re-Entry Trajectory, Applied Sciences, 14 (2024): 10585.

[4] I. Malmir, Novel closed-loop controllers for fractional linear quadratic time-varying systems, Numerical Algebra, Control and Optimization, 14, no. 2 (2024): 366–403.

[5] S. Kumar, A. Kumar Awasthi, S. Kumar Mishra, H.Chandra Yadav, S. Lal, An error estimation of absolutely continuous signals and solution of Abel's integral equation using the first kind pseudo-Chebyshev wavelet technique, Franklin Open, 10 (2025): 100205.

[6] M. Muhafzan, A. Nazra, A. I. Baqi, Z. Zulakmal, On tracking control problem for fractional order descriptor systems, CYBERNETICS AND PHYSICS, 11 (2022): 87–93.

[7] C. Han, Y. Chen, D.-Y. Liu, and D. Boutat, Numerical analysis of viscoelastic rotating beam with variable fractional order model using shifted Bernstein–Legendre polynomial collocation algorithm, Fractal and Fractional 5, no. 1 (2021): 8.

[8] A. K. Hameed, M. M. Mustafa, Numerical Solution of Linear Fractional Differential Equation with Delay Through Finite Difference Method, Iraqi Journal of Science, (2022): 1232-1239.

[9] J. Sheng-Hong Tsai, Y.-T. Liao, F. Ebrahimzadeh, S.-Y. Lai, T.-J. Su, S.-M. Guo, L.-S. Shieh, and T.-J. Tsai, A new PI optimal linear quadratic state-estimate tracker for continuous-time non-square non-minimum phase systems, International Journal of Systems Science, 48 (2017): 1438–1459.

[10] M. Pal, S. Mondal, S. Chowdhury, A. Mondal, A. Mukherjee, Core Power Control of a Pressurized Water Reactor Using Event Triggered PID Control, In 2025 IEEE 1st International Conference on Smart and Sustainable Developments in Electrical Engineering (SSDEE), (2025): 1–6.

[11] I. Malmir, Novel closed-loop controllers for fractional nonlinear quadratic systems, Mathematical Modelling and Control, 3 (2023): 345–354.

[12] V. Ayala, M. Torreblanca, W. Valdivia, Toward Applications of Linear Control Systems on the Real World and Theoretical Challenges, Symmetry, 13, no. 2 (2021): 167.

[13] Y. Tian, Developing Computational Models of Basal Ganglia and Thalamic Nuclei Impacted by Deep Brain Stimulation (DBS), Implications for Model-based Closed-loop DBS Control, PhD diss., University of Toronto, Canada, 2024.

[14] M. Muhafzan, A. Nazra, A. I. Baqi, Z. Zulakmal, On the fractional tracking control and its application in diabetes treatment, In AIP Conference Proceedings, vol. 2920, no. 1, AIP Publishing, (2024).

[15] I. Malmir, Suboptimal control law for a multi fractional high order linear quadratic regulator system in the presence of disturbance, Results in Control and Optimization, (2023), 100251.

[16] I. Malmir, A novel wavelet-based optimal linear quadratic tracker for time-varying systems with multiple delays, Statistics, Optimization & Information Computing, 9, no. 2 (2021): 418–434.

[17] I. Malmir, A general framework for optimal control of fractional nonlinear delay systems by wavelets, Statistics, Optimization & Information Computing, 8, no. 4 (2020): 858–875.

[18] J. Moon, The Pontryagin type maximum principle for Caputo fractional optimal control problems with terminal and running state constraints. AIMS Mathematics, 10 (2025): 884–920.

[19] J.-P. Huang, H.-C. Zhou, Infinite Horizon Linear Quadratic Optimal Control Problems for Singular Volterra Integral Equations, SIAM Journal on Control and Optimization, 63 (2025): 57–85.

[20] G. Dewangan, A. Singh, A. Kanaujiya, Generalized distributed-order fractional optimal control problem using Laguerre wavelet method, Journal of Mathematical Modeling, (2025): 707–725.

[21] A. Singh, A Kanaujiya., J. Mohapatra, Chelyshkov wavelet method for solving multidimensional variable order fractional optimal control problem, Journal of Applied Mathematical and Computing, 70 (2024): 3135–3160.

[22] S. Karami, M. H. Heydari, D. Baleanu, M. Shasadeghi, Generalized Hat Functions for Fractional Delay Optimal Control Problems With $\psi$-Caputo Fractional Derivative, Mathematical Methods in the Applied Sciences, (2025).

[23] I. Malmir, An efficient method for a variety of fractional time-delay optimal control problems with fractional performance indices, International Journal of Dynamics and Control, 11, no. 2 (2023): 2886–2910.

[24] I. Malmir, New pure multi-order fractional optimal control problems with constraints: QP and LP methods, ISA transactions, 153 (2024): 155–190.

[25] I. Malmir, Pure Fractional Optimal Control of Partial Differential Equations: Nonlinear, Delay and Two-Dimensional PDEs, Pan-American Journal of Mathematics, 4 (2025): 9.

[26] R. S. Hilscher, V. M. Zeidan, Transformation preserving controllability for nonlinear optimal control problems with joint boundary conditions, ESAIM: Control, Optimisation and Calculus of Variations, 27 (2021): 75.

[27] I. Malmir, S. H. Sadati, Transforming linear time-varying optimal control problems with quadratic criteria into quadratic programming ones via wavelets, Journal of Applied Analysis, 26 (2020): 131–151.

[28] S. Lal, D. K. Singh, Moduli of continuity of functions in Holder's class and solution of Rayleigh, Van der Pol and Duffing equations by sixth-kind Chebyshev wavelets, Indian Journal of Pure and Applied Mathematics, (2025).

[29] T. Cuchta, D. Poulsen, N. Wintz, Linear Quadratic Tracking with Continuous Conformable Derivatives, European Journal of Control, 72 (2023): 100808.

[30] M. L. Nagurka and S.-K. Wang, A Chebyshev-based state representation for linear quadratic optimal control, Journal of Dynamic Systems, Measurement, and Control, 115 (1993): 1–6.

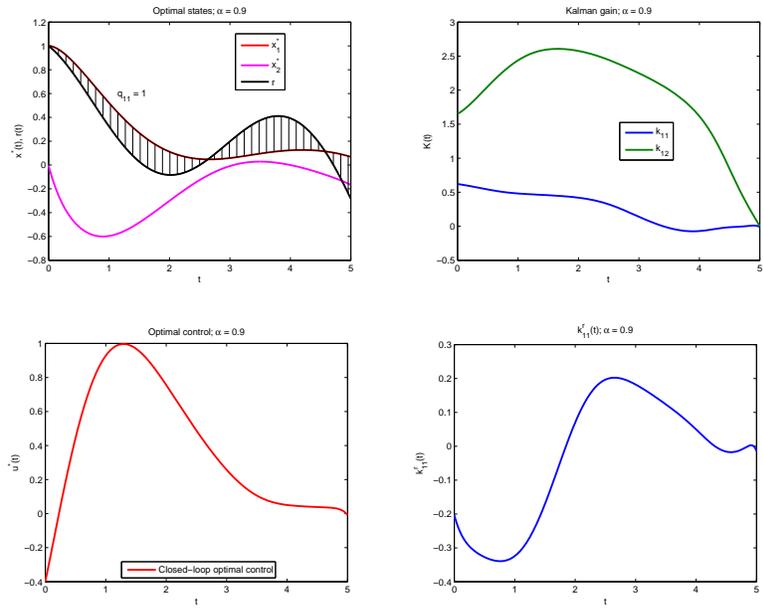

Figure 1: Optimal solutions for Example 1, $q_{11} = 1$ and $\alpha = 0.9$

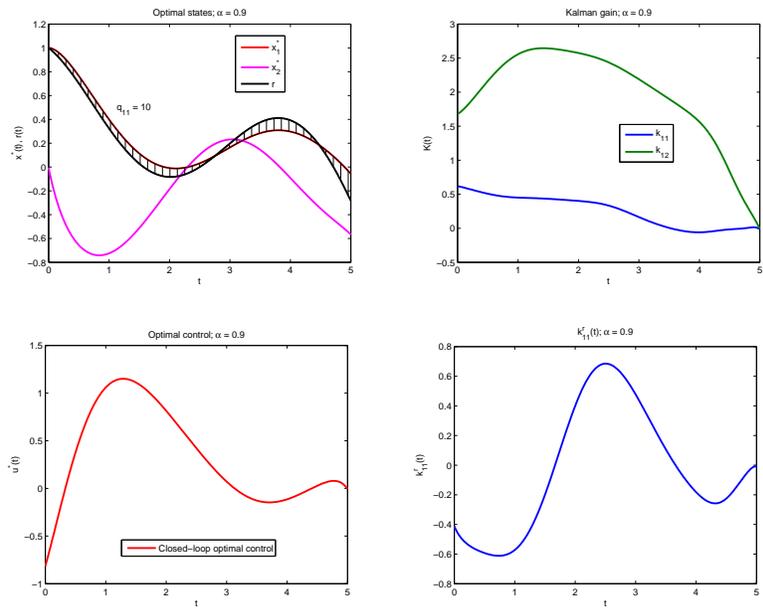

Figure 2: Optimal solutions for Example 1, $q_{11} = 10$ and $\alpha = 0.9$

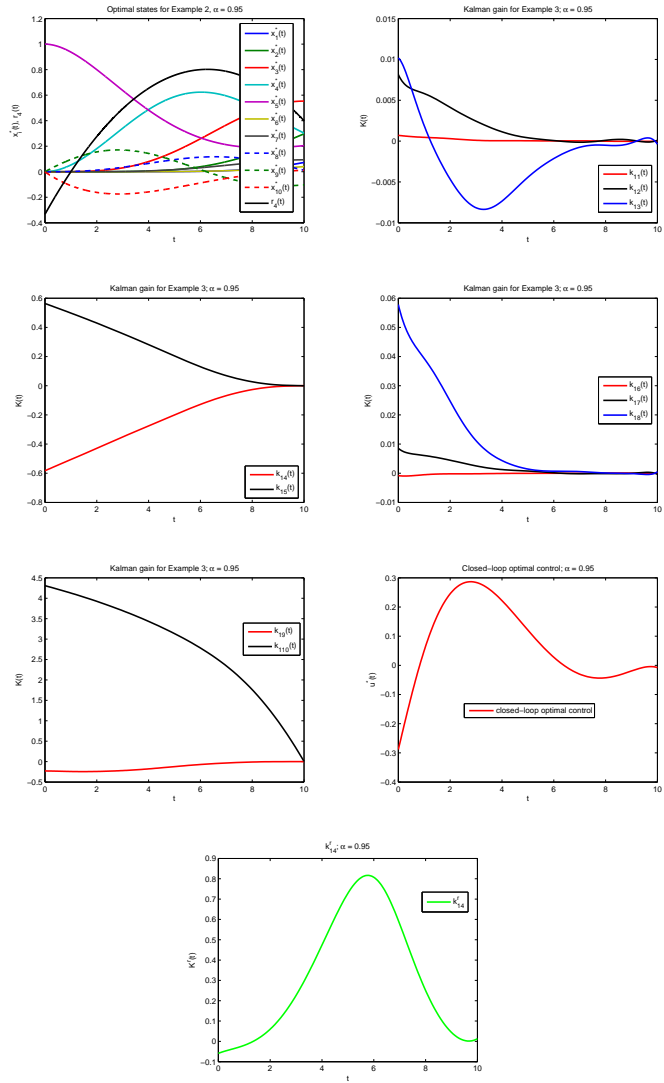

Figure 3: Optimal solutions for Example 2 and $\alpha = 0.95$